\documentclass[a4paper,12pt]{article}

\usepackage[latin1]{inputenc}
\usepackage[english]{babel}
\usepackage{amsmath}
\usepackage{amsfonts}
\usepackage{amssymb}
\usepackage{amsthm}
\usepackage{graphicx}

\title{Heegaard splittings of Seifert manifolds without an orientable base space}
\author{Simone Penzavalle\footnote{Affiliation and mailing address: \newline\indent Scuola Normale Superiore \newline\indent Piazza dei Cavalieri 7 \newline\indent 56126 Pisa \newline\indent Italy}}
\date{}

\begin{document}

\newtheorem{thm}{Theorem}[section]
\newtheorem{defin}[thm]{Definition}
\newtheorem{lem}[thm]{Lemma}
\newtheorem{prop}[thm]{Proposition}
\newtheorem{es}[thm]{Example}
\newtheorem{remark}[thm]{Remark}

\newcommand{\esseuno}{S\hspace{0.5mm}^1}
\newcommand{\essedue}{\mathbb{S}^2}
\newcommand{\didue}{\mathbb{D}^2}
\newcommand{\I}{{\rm \bf I}}
\newcommand{\fti}{F \tilde{\times} \I}
\newcommand{\pidue}{\mathbb{RP}^2}
\newcommand{\linea}{\vspace{1em}}

\newenvironment{rem}{\begin{remark} \upshape}{\end{remark}}
\newenvironment{defi}{\begin{defin} \upshape}{\end{defin}}

\maketitle

\begin{small}
\noindent\textbf{Abstract.} We make use of Heath's, Moriah's and Schultens's results to prove that irreducible Heegaard splittings of orientable Seifert manifolds with nonorientable base space are either vertical or horizontal.

\begin{center}

\end{center}
\end{small}

\section*{Introduction}

The purpose of this paper is to complete the characterization of irreducible Heegaard splittings of orientable Seifert fibered spaces which was given by Moriah and Schultens in \cite{moriahschultens-seifert} for the totally orientable case to the case of orientable Seifert fibered spaces whose fibration is not orientable.

The main Theorem of the paper, which is analogous to \cite[Theorem 0.1]{moriahschultens-seifert}, is the following.

\linea

\noindent\textbf{Theorem \ref{main}.} \emph{Let $ M $ be an orientable Seifert fibered space with nonorientable base space. Then every irreducible Heegaard splitting of $ M $ is either vertical or horizontal.}

\linea

In order to prove Theorem \ref{main} we make use of the results Schultens proved in \cite{schultens-graph}, after equipping a Seifert fibered space with boundary with a totally orientable generalized graph structure.

The paper is organized as follows.

Section 1 is devoted to the definition of Heegaard splittings of a compact {\rm 3}-manifold $ M $ and to the recollection of known properties of such structures. In particular, we recall the notion of \emph{weak reduction} of a pair $ (M, \Sigma) $, where $ \Sigma $ is a Heegaard surface in $ M $, and that of \emph{amalgamation} of two such pairs. We also recall Scharlemann's and Thompson's characterization of Heegaard splittings of handlebodies and of products $ F \times \I $, with $ F $ a closed orientable surface, and Heath's result about Heegaard splittings of twisted $ \I $-bundles.

Section 2 deals with orientable Seifert fibered spaces which do not admit a totally orientable Seifert structure. A natural Heegaard splitting structure on such manifolds is that of a \emph{vertical} splitting, which is defined as in the totally orientable case. But there are manifolds where the structure of a \emph{horizontal} splitting can arise, and the two structures are often distinguished by the genus of the Heegaard surface. The definition of horizontal splitting is analogous to that given in the totally orientable case, taking care of the different properties of horizontal surfaces.

We first focus on strongly irreducible splittings of manifolds with nonempty boundary. We equip such a manifold with a totally orientable generalized graph structure, and then use Schultens's \cite[Theorem 1.1]{schultens-graph} to find a suitable isotopy of the Heegaard surface. We obtain the following Theorem.

\linea

\noindent\textbf{Theorem \ref{strirrboundary}.} \emph{Any strongly irreducible Heegaard splitting of a Seifert fibered space with nonorientable base and nonempty boundary is vertical.}

\linea

Strongly irreducible Heegaard splittings of manifolds with boundary are known to be the ``building blocks'' of weakly reducible splittings. More precisely, any weakly reducible splitting of any manifold $ M $ (with boundary or not) can be obtained by amalgamation of strongly irreducible splittings of manifolds which can be obtained by cutting $ M $ along suitable incompressible surfaces. We apply this idea and prove the following Theorem.

\linea

\noindent\textbf{Theorem \ref{irrweakredsplit}.} \emph{Any weakly reducible Heegaard splitting of a Seifert fibered space with nonorientable base is vertical.}

\linea

The only case left out is that of strongly irreducible splittings of closed Seifert manifolds. Here we follow very closely Moriah's and Schultens's paper \cite{moriahschultens-seifert}. First, using strong irreducibility one shows that it is possible to push some fiber $ f $ into the Heegaard surface $ \Sigma $. Then one looks at what $ \Sigma $ looks like in the complement of a fibered neighbourhood of $ f $ in $ M $. It turns out that if it is incompressible then either the given splitting is horizontal or the manifold $ M $ is a lens space, and that if it is compressible then the given splitting is vertical. As Heegaard splittings of lens spaces are already characterized (\cite{bonahonotal-heegaardlens}, \cite{moriahschultens-seifert}), we obtain Theorem \ref{main}.

I wish to thank Luisa Paoluzzi for many helpful discussions. I also wish to thank Carlo Petronio for suggesting me to work on this subject.

\section{Preliminaries}

This section is devoted to the description of Heegaard splittings and their properties. In particular, Lemma \ref{weaklyred}, due to Scharlemann, is of importance. We define the operations of \emph{weak reduction} of a pair $ (M, \Sigma) $, where $ M $ is a compact {\rm 3}-manifold and $ \Sigma $ a Heegaard surface in $ M $, and of \emph{amalgamation} of two such pairs, and finally recall Scharlemann and Thompson's results on irreducible splittings of handlebodies and products $ (closed $ $ orientable $ $ surface) \times \I $.

\subsection{Compression bodies and Heegaard splittings}

A \emph{compression body} $ W $ is a {\rm 3}-manifold constructed
by gluing {\rm 2}-handles to $ S \times \I $, where $ S $ is a
closed orientable connected surface, along a collection of
disjoint simple closed curves contained in $ S \times \{ 0 \} $,
and capping off any sphere component originating from $ S \times \{
0 \} $ with a {\rm 3}-handle. The connected component $ S \times
\{ 1 \} $ of $ \partial W $ is denoted by $ \partial_+ W $, while
$ \partial_- W $ denotes $ \partial W \setminus \partial_+ W $. A
compression body $ W = \partial_+ W \times \I $ is called
\emph{trivial}, while a compression body $ W $ such that $
\partial_- W = \emptyset $ is called a \emph{handlebody}. The
\emph{genus} of $ W $ is the genus of $ \partial_+ W $.

There is a dual method for constructing compression bodies. One
can start with a (perhaps empty) closed orientable non necessarily connected surface $ S $ with no sphere components and a ball $ B $ and glue $ k $ 1-handles on $ (S \times \{ 1 \}) \cup \partial B \subset (S
\times \I) \cup B $ in such a way that the resulting manifold is
connected. Then the resulting manifold is a compression body $ W $
with $ \partial_- W = S \times \{ 0 \} $. The cocores of the 1-handles are called
\emph{meridian disks}. A collection $ \Upsilon $ of meridian disks for a compression body $ W $ is called \emph{complete} if $ W \setminus \Upsilon $ is either a ball or a trivial compression body.

\begin{rem} \label{cprbdy}
Let $ W $ be a compression body, and $ \Delta $ a collection of properly embedded disks in $ W $. Then the {\rm 3}-manifold $ W' $ obtained from $ W $ by cutting along $ \Delta $ is a union of compression bodies. If in addition $ \partial \Delta \subset \partial_+ W $, then $ \partial_- W' = \partial_- W $.
\end{rem}

A properly embedded {\rm 1}-complex $ Q $ such that $ W $ collapses to $ Q \cup \partial_- W $ is called a \emph{spine} of $ W $. A simple closed curve contained in a spine of $ W $ is called a \emph{core} of $ W $.

Spines are related to meridian disks in the following way. Given a meridian disk $ D $, the core of the 1-handle whose cocore is $ D $ is an edge of a spine $ Q $, namely the spine formed by the cores of the 1-handles used to build $ W $. Vice-versa, given a spine $ Q $, as $ W $ is a regular neighbourhood of $ Q \cup \partial_- W $, it can be constructed by gluing 1-handles to $ \partial_- W \times \I $ using the edges of $ Q $ as cores, and the cocores of these handles are meridian disks.

\begin{defi}
Given a {\rm 3}-manifold $ M $, a closed connected properly
embedded orientable surface $ \Sigma $ is a \emph{Heegaard
surface} (or a \emph{splitting surface}) for $ M $ if $ M
\setminus \Sigma $ consists of two compression bodies $ W_1 $, $
W_2 $ with $ \partial_+ W_1 = \partial_+ W_2 = \Sigma $. The pair
$ ( W_1 , W_2 ) $ is also called a \emph{Heegaard splitting} of $
M $.

The \emph{genus} of $ ( W_1 , W_2 ) $ is the genus of $ \Sigma $.
\end{defi}

In the following, Heegaard surfaces will always be considered up
to isotopy.

\begin{defi}
Given a {\rm 3}-manifold $ M $ and a Heegaard surface $ \Sigma $
for $ M $, the surface $ \Sigma' $ obtained by taking the
connected sum of pairs $ ( M , \Sigma ) \# ( S^3 , T ) $, where $
T $ is a torus splitting $ S^3 $ into two solid tori, is called
the \emph{stabilization} of $ \Sigma $. The Heegaard surface constructed this way is \emph{stabilized}.
\end{defi}

An equivalent condition for a Heegaard splitting to be stabilized
is the existence of properly embedded disks $ D_1 \subset W_1 $
and $ D_2 \subset W_2 $ with $ \partial D_1 \cap \partial D_2 = \{
p \} \subset \Sigma $.

\begin{defi}
A Heegaard splitting is called \emph{reducible} if there are
properly embedded essential disks $ D_1 \subset W_1 $ and $ D_2
\subset W_2 $ with $ \partial D_1 = \partial D_2 $. Otherwise the Heegaard splitting is \emph{irreducible}.

A Heegaard splitting is called \emph{weakly reducible} if there
are properly embedded essential disks $ D_1 \subset W_1 $ and $
D_2 \subset W_2 $ with $ \partial D_1 \cap \partial D_2 =
\emptyset $. Otherwise the Heegaard splitting is \emph{strongly irreducible}.
\end{defi}

The following Theorem collects useful well-known properties of Heegaard splittings. Statements {\rm 1} and {\rm 2} are straightforward (see e.g. \cite{scharlemann-heegaard}), {\rm 3} follows from Waldhausen's results on Heegaard splittings of $ S^3 $ (\cite{waldhausen-heegaardsphere}), {\rm 4} is due to Haken (\cite{haken-results3mflds}), {\rm 5} is just a generalization of {\rm 4}, and {\rm 6} is an easy application of Morse theory (see e.g. \cite{milnor-morse}).

\begin{thm} \label{propsplit}
$ $
\begin{enumerate}
   \item A reducible Heegaard splitting is weakly reducible.
   \item A stabilized Heegaard splitting different from the genus
   {\rm 1} splitting of $ S^3 $ is reducible.
   \item A reducible Heegaard splitting of an irreducible {\rm
   3}-manifold is stabilized.
   \item \label{ridrid} Any Heegaard splitting of a reducible {\rm 3}-manifold is
   reducible.
   \item \label{boundridweakrid} Any Heegaard splitting of a $
   \partial $-reducible {\rm 3}-manifold is either weakly
   reducible or one compression body is trivial.
   \item To every Heegaard splitting of a {\rm 3}-manifold $ M $
   is associated a Morse function $ h : M \rightarrow [0,1] $ with
   distinct critical values
   $$ 0 < a_1 < ... < a_k < b_1 < ... < b_l < 1 $$
   such that:
   \begin{itemize}
      \item $ h^{-1}(a_i) $ is a critical point of index $ 1 $ for
      all $ i = 1, ..., k $ and $ h^{-1}(b_j) $ is a critical
      point of index $ 2 $ for all $ j = 1, ..., l $;
      \item $ h^{-1}(0) = \partial_- W_1 $ if this is nonempty, or an index {\rm 0} critical
		point if $ \partial_- W_1 = \emptyset $; and $ h^{-1}(1) = \partial_- W_2 $ if this is
		nonempty, or an index {\rm 3} critical point if $ \partial_- W_2 = \emptyset $;
      \item $ \Sigma = h^{-1}(r) $ with $ a_k < r < b_1 $.
   \end{itemize}
\end{enumerate}
\end{thm}

The following result can be found in \cite[Lemma 3.10]{scharlemann-heegaard}.

\begin{lem}[Scharlemann] \label{weaklyred}
Let $ \Sigma $ be a strongly irreducible Heegaard surface in a
{\rm 3}-manifold $ M $, and $ D $ a disk in $ M $ with $ \partial
D \subset \Sigma $ with transverse intersection. Then $ D $ can be
isotoped to a disk in $ W_1 $ or $ W_2 $.
\end{lem}

\subsection{Amalgamation and weak reductions of Heegaard splittings}

Throughout the whole paper, for $ X $ a submanifold or a complex in a manifold $ M $, by $ \eta(X) $ we will denote a regular neighbourhood of $ X $ in $ M $.

Let $ M $, $ N $ be {\rm 3}-manifolds, $ B_M $ a connected component of $ \partial M $ and $ B_N $ a connected component of $ \partial N $, and $ h : B_M \rightarrow B_N $ a homeomorphism. Let $ (W^M_1, W^M_2) $ be a Heegaard splitting of $ M $ with $ B_M \subset \partial_- W^M_1 $ and $ (W^N_1, W^N_2) $ a Heegaard splitting of $ N $ with $ B_N \subset \partial_- W^N_1 $. We want to exploit the product structure near $ B_M $ and $ B_M $ to produce a Heegaard splitting of the manifold $ M \cup_h N $.

Set $ B'_M = \partial_- W^M_1 \setminus B_M $ and $ B'_N = \partial_- W^N_1 \setminus B_N $, so that $ W^M_1 $ is constructed by attaching {\rm 1}-handles to $ \bar{\eta}(B_M \cup B'_M) $ and $ W^N_1 $ is constructed attaching {\rm 1}-handles to $ \bar{\eta}(B_N \cup B'_N) $. Let $ f : \bar{\eta}(B_M)
\rightarrow B_M \times \I $ and $ g : \bar{\eta}(B_N) \rightarrow
B_N \times \I $ be homeomorphisms. Consider the following
equivalence relation on $ M \sqcup N $:
\begin{itemize}
   \item if $ x, y \in \eta(B_M) $, then $ x \sim y $ iff $ p_M(f(x))
   = p_M(f(y)) $;
   \item if $ x, y \in \eta(B_N) $, then $ x \sim y $ iff $ p_N(g(x))
   = p_N(g(y)) $;
   \item if $ x \in B_M $ and $ y \in B_N $, then $ x \sim y $ iff $ y
   = h(x) $;
\end{itemize}
where $ p_M $ and $ p_N $ denote the projection onto the first coordinate in $ B_M \times \I $ and $ B_N \times \I $ respectively.

Up to isotopies, we can assume that the attaching disks for $ W^M_1 $'s and $ W^N_1 $'s {\rm 1}-handles have disjoint equivalence classes. Set
$$ M \cup_h N = \frac{M \cup N}{\sim}, \qquad W_1 = \frac{W^M_1 \cup
W^N_2}{\sim}, \qquad W_2 = \frac{W^M_2 \cup W^N_1}{\sim}. $$

Then $ W_1 $ is obtained from $ W^N_2 \cup \bar{\eta}(B'_M) $ by
connecting $ \partial_+ W^N_2 $ to $ f^{-1}(B'_M \times \{ 1 \}) $
via {\rm 1}-handles, and so is a compression body. Analogously, $
W_2 $ is obtained from $ W^N_1 \cup \bar{\eta}(B'_N) $ by connecting
$ \partial_+ W^N_1 $ to $ g^{-1}(B'_N \times \{ 1 \}) $ via {\rm
1}-handles. Hence $ (W_1, W_2) $ is a Heegaard splitting of $ M
\cup_h N $.

\begin{defi}
The Heegaard splitting $ (W_1, W_2) $ is called the \emph{amalgamation} of $ (W^M_1, W^M_2) $ and $ (W^N_1, W^N_2) $ along $ B_M, B_N $ via $ h
$.
\end{defi}

We now reformulate a result which comes from \cite{scharlthompson-weakreduction}.

\begin{thm}[Scharlemann, Thompson] \label{weakreduction}
Let $ \Sigma $ be an irreducible Heegaard surface in a {\rm 3}-manifold $ M $. Then, setting $ F_0 = \partial_- W_1 $ and $ F_m = \partial_- W_2 $, there exist properly embedded disjoint connected surfaces $ \Sigma_1, ..., \Sigma_m $ and $ F_1, ..., F_{m - 1} $ such that:
\begin{itemize}
   \item $ F_i $ is incompressible for $ i = 1, ..., m - 1 $;
   \item all $ F_i $'s belong to the same homology class in $ H_2(M) $: let $ M_i $ be the cobordism between $ F_{i - 1} $ and $ F_i $ for $ i = 1, ..., m $;
   \item $ \Sigma_i $ is a strongly irreducible Heegaard splitting of $ M_i $;
   \item $ \Sigma $ is obtained by amalgamation of the $ \Sigma_i $'s along the $ F_i $'s.
\end{itemize}
Such a decomposition will be called a \emph{weak reduction} of $ \Sigma $. If $ \Sigma $ itself is strongly irreducible (e.g., if $ M $ is not Haken) then one can take $ m = 1 $.
\end{thm}

The following Theorems are respectively \cite[Lemma 2.7]{scharlthompson-fxi}, \cite[Theorem 2.11]{scharlthompson-fxi} and \cite[Corollary 1.6]{heath-classificationheegaard}.

\begin{thm}[Scharlemann, Thompson] \label{handlebody}
Any Heegaard splitting of a handlebody $ H $ is \emph{standard},
i.e. it is obtained by repeated stabilization of the
\emph{trivial} splitting in which $ \Sigma $ is parallel to $
\partial H $.
\end{thm}

\begin{thm}[Scharlemann, Thompson] \label{fxi}
Let $ F \neq \essedue $ be a closed orientable surface. Let
$ \Sigma $ be an irreducible Heegaard surface for $ F \times \I $.
Then either:
\begin{enumerate}
   \item $ \Sigma = F \times \{ 1/2 \}$, or
	\item the genuine compression body $ W_2 $ admits a spine consisting of a single vertical
	arc $ \{ p \} \times \I $.
\end{enumerate}
\end{thm}

\begin{thm}[Heath] \label{heath}
Let $ F $ be a closed nonorientable surface. Then every Heegaard splitting of $ \fti $ is obtained by repeated stabilization of the \emph{trivial} splitting obtained by taking a boundary parallel surface and adding a vertical {\rm 1}-handle.
\end{thm}

\section{Seifert manifolds}

We begin by recalling a few known facts. The following Proposition follows from \cite[Lemma 1.7]{jaco-3mfld}.

\begin{prop}[Jaco] \label{excluded}
The only Seifert manifolds which are not irreducible are $ \essedue \times \esseuno $ and $ \pidue \tilde{\times} \esseuno \cong \mathbb{RP}^3 \# \mathbb{RP}^3 $.
\end{prop}
\begin{rem} \label{onlylarge}
In the whole of the present work we will assume our manifold to be irreducible and to admit an incompressible vertical torus. This is not a gap in the characterization because:
\begin{itemize}
   \item any Heegaard splitting of $ \mathbb{RP}^3 \# \mathbb{RP}^3 $ is reducible (see \ref{ridrid} of Theorem \ref{propsplit}) and can be obtained by Heegaard splittings of $ \mathbb{RP}^3 $ by a standard construction (see e.g. \cite{scharlemann-heegaard}), and \cite{bonahonotal-heegaardlens} characterizes Heegaard splittings of lens spaces, in particular those of $ L_{2,1} \cong \mathbb{RP}^3 $;
   \item the only orientable manifolds which admit a nonorientable Seifert structure without incompressible vertical tori are those with base orbifold $ \pidue $ and at most one exceptional fiber, and these (see the Remark following \cite[Theorem 6.2.2]{orlik-seifert}) are either $ \mathbb{RP}^3 \# \mathbb{RP}^3 $ or the \emph{prism manifolds} $ M(\essedue, e, (2,\beta_1), (2,\beta_2), (\alpha_3,\beta_3)) $. Thus, up to changing the Seifert structure, irreducible Heegaard splittings are already characterized for these manifolds (\cite{bonahonotal-heegaardlens},\cite{moriahschultens-seifert}).
\end{itemize}
\end{rem}

We now state the results proved in \cite{schultens-graph}, in order
to use them to describe Heegaard splittings of Seifert manifolds
with nonorientable base.

\begin{defi}
A surface $ \Sigma $ in a Seifert fibered manifold $ M $ is \emph{vertical} if it is saturated with regard to the Seifert fibration; it is \emph{horizontal} if it intersects every fiber transversely, and it is \emph{pseudohorizontal} if there is a fiber $ f \subset \Sigma $ such that $ \Sigma \cap (M \setminus \eta(f)) $ is horizontal in $ M \setminus \eta(f) $. A surface $ \Sigma $ is called \emph{pseudovertical} if it can be obtained by a collection $ \mathcal{A} $ of properly embedded vertical surfaces via {\rm 1}-surgery along a collection $ \Gamma $ of arcs with endpoints on $ \mathcal{A} $ which projects to an embedded collection of arcs in $ S $.
\end{defi}

Given some pseudovertical surface $ \Sigma $, by general position we can avoid intersections between singular fibers and $ \Sigma $'s graph $ \Gamma $. Setting $ \mathcal F = \{ $\emph{singular fibers}$ \} $ if this is nonempty, and $ \mathcal F = \{ $\emph{one regular fiber}$ \} $ otherwise, as $ M $ admits a section out of a neighbourhood of $ \mathcal F $, we can identify $ S^* = S \setminus \pi(\eta(\mathcal F)) $ with its image via this section. Then, as each arc is contractible, we can isotope $ \Gamma $ until it lies in $ S^* $.

As the following Lemma shows, it is very easy to produce new representations of a given pseudovertical surface that correspond to isotopies.

\begin{lem} \label{lemswap}
Let Figure \ref{swap} (left) represent $ S^* $, where solid curves represent its intersection with $ \mathcal A $ and dashed arcs represent $ \Gamma $.

\begin{figure}[!htb]
\begin{center}
\input{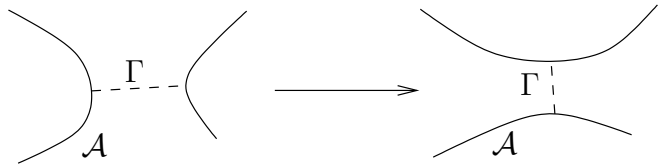}
\caption{Swapping along an arc.} \label{swap}
\end{center}
\end{figure}

Then changing the configuration to that of Figure \ref{swap} (right) corresponds to an isotopy of the corresponding pseudovertical surface.
\end{lem}
\begin{proof}
The required isotopy consists in pushing the chosen arc along the vertical direction until the whole fiber has been traversed. Once more we remark that this is possible because the fibration is trivial over $ \Gamma $.
\end{proof}

Hence we have a move that we can perform on the system of curves and arcs $ (\mathcal{A} \cup \Gamma) \cap S^* $ that changes the pseudovertical surface via an isotopy. Such a move will be called a \emph{swap}. Note that if $ \alpha' $ is the arc obtained performing a swap along some arc $ \alpha $, when we perform a swap along $ \alpha' $ we get back the same arc $ \alpha $ and the same vertical surfaces.

\begin{rem} \label{swapmoebius}
Let $ M $ be a M\"oebius band contained in the base $ S $, and suppose $ \mathcal A $ contains $ \pi^{-1}(\partial M) $. Suppose further that one of the arcs in $ \Gamma $ (call it $ \gamma $) is entirely contained in $ M $. Then, either $ \gamma $ cobounds a disk together with a subarc of $ \partial M $, or $ \gamma $ is essential in $ M $. In the latter case, swapping along $ \gamma $ produces an inessential vertical torus (because cutting $ M $ along $ \gamma $ gives a disk) and an arc $ \gamma' $ outside it. As the topology of $ S $ has not changed, we conclude that swapping along an arc which is essential in some cross-cap of $ S $ has merely the effect of pushing the cross-cap through $ \mathcal A $.
\end{rem}

We will now recall the definition of \emph{vertical} Heegaard splittings of a Seifert fibered manifold.

Let $ M $ be a Seifert fibered manifold with genus-$ g $ base $ S $, exceptional fibers $ \mathcal F = \{ f_1, ..., f_m \} $ and boundary components $ \mathcal B = \{ B_1, ..., B_n \} $; and let $ \mathcal E = \{ e_1, ..., e_m \} $ be the set of exceptional points on $ S $ and $ \{ b_1, ..., b_n \} $ the set of boundary components of $ S $. Finally, let $ \sigma_i $, for $ i = 1, ..., n $, be a parallel copy of $ b_i $ such that the annulus cobounded by $ \sigma_i $ and $ b_i $ contains no exceptional points, and let $ \Sigma_i $ be the vertical torus over $ \sigma_i $. Then there are three possibilities:
\begin{itemize}
   \item if $ n > 0 $, i.e. if $ \partial M \neq \emptyset $, fix a boundary component $ \overline B $ and partition $ \mathcal B \setminus \{ \overline B \} $ into subsets $ \mathcal B_1 $ and $ \mathcal B_2 $ and $ \mathcal E $ into subsets $ \mathcal E_1 $ and $ \mathcal E_2 $;
   \item if $ n = 0 $ and $ m > 0 $, then fix an exceptional fiber $ \bar f $ and partition $ \mathcal E \setminus \{ \bar f \} $ into subsets $ \mathcal E_1 $ and $ \mathcal E_2 $, and set $ \mathcal B_1 = \mathcal B_2 = \emptyset $;
   \item if $ n = m = 0 $ then set $ \mathcal B_1 = \mathcal B_2 = \mathcal E_1 = \mathcal E_2 = \emptyset $.
\end{itemize}

Now fix a regular base point $ p \in S $ and choose a properly embedded graph $ \Gamma \subset S \setminus \{ \sigma_i | B_i \in \mathcal B_1 \} $ having $ \{ p \} \cup \mathcal E_1 $ as vertex set and such that $ S \setminus ( \Gamma \cup \{ \sigma_i | B_i \in \mathcal B_1 \} ) $ is a regular neighbourhood of $ \partial S \cup (\mathcal E \setminus \mathcal E_1) $ or of a point if this set is empty.

Note that the Seifert fibration admits a section when restricted over $ \Gamma $; we now identify $ \Gamma $ with its image in $ M $ via this section. If both $ n, m = 0 $ add to $ \Gamma $ the regular fiber over $ p $. Let $ \Sigma $ be the surface obtained via {\rm 1}-surgery on $ \{ \Sigma_i | B_i \in \mathcal B_1 \} $ along $ \Gamma \cup \mathcal E_1 $.

It follows from the very definition that the surface $ \Sigma $ is a pseudovertical surface. It is not difficult to prove, using the swap move, that $ \Sigma $ is indeed a Heegaard surface.

\begin{defi} \label{defvertical}
A Heegaard splitting constructed as above is called \emph{vertical}.
\end{defi}

The swap move allows also to prove the following Lemma, which is already known in literature (\cite[Proposition 2.10]{schultens-fpers1}).

\begin{lem}[Schultens] \label{amalgvert}
Let $ M_1 $ and $ M_2 $ be Seifert manifolds with surfaces $ \Sigma_1 $ and $ \Sigma_2 $ both giving vertical Heegaard splittings. Then the Heegaard splitting obtained by amalgamation of the $ \Sigma_i $'s via a fiberwise homeomorphism between some boundary components of the $ M_i $'s is vertical.
\end{lem}

The following Proposition is very powerful, and will be used later.

\begin{prop} \label{equivvert}
An irreducible Heegaard splitting of a Seifert fibered manifold with base orbifold $ S $ is vertical if and only if the corresponding Heegaard surface is pseudovertical.
\end{prop}

\begin{rem} \label{comprseifert}
If a compression body $ W $ admits a Seifert fibration, then either $ W = (annulus) \times \esseuno $ or $ W $ is a solid torus.
\end{rem}

\begin{proof}[Proof of Proposition \ref{equivvert}]
We only have to prove the ``if'' part of the Proposition. First of all, note that given a separating pseudovertical surface $ \Sigma $, the corresponding family of vertical surfaces $ \mathcal A $ is also separating. This is because it is obtained from $ \Sigma $ by repeated compression of {\rm 1}-handles, and at each stage the homology class of the surface does not change. In particular, compressing a separating surface yields a (perhaps disconnected) separating surface.

Recall that thanks to general position we ensured that singular fibers are disjoint from $ \Gamma $, and observe that no isotopy can make one of them lie in $ \mathcal A $, or we would lose (proper) embeddingness.

Now perform swaps so that $ \Gamma $ lies entirely on the same side of $ \mathcal A $, and cut $ M $ along $ \mathcal A $. What we get is a connected component $ S_\Gamma $ containing $ \Gamma $, and a (perhaps disconnected) part $ S' $ not containing $ \Gamma $. Now, $ \pi^{-1}(S') $ is Seifert fibered and is obtained by compressing a compression body along the cocores of its {\rm 1}-handles, hence it is a union of $ (annuli) \times \esseuno $ and solid tori thanks to Remark \ref{cprbdy} and Lemma \ref{comprseifert}; hence $ S' $ consists of annuli without exceptional points and disks with at most one exceptional point each.

Now cut $ S_\Gamma $ along $ \Gamma $. To see which pieces one gets, we can swap along $ \Gamma $ and look at the connected components not containing the new set of arcs. This again corresponds to perform compressions of the {\rm 1}-handles of a compression body. Hence, these pieces are again annuli without exceptional points and disks with at most one exceptional points each.

Hence cutting along $ \Gamma $ kills the whole genus of $ S $ and does not leave any two components of $ \partial S \cup \mathcal F $ in the same complementary connected component, thus the splitting is vertical.
\end{proof}

\subsection{Graph structures} \label{structure}

We will now make use of some definitions and results of \cite{schultens-graph}.

\begin{defi}
Let $ M $ be a {\rm 3}-manifold. A \emph{generalized graph structure} on $ M $ is a finite graph $ \Gamma $, with vertex set $ V $ and edge set $ E $, such that:
\begin{itemize}
   \item to each $ v \in V $ we can associate a {\rm 3}-manifold $ M_v $
   which is either Seifert fibered or homeomorphic to $ S_v \times \I $,
   with $ S_v $ a compact surface;
   \item to each $ e \in E $ we can associate a {\rm 3}-manifold $ M_e $
   homeomorphic either to $ T^2 \times \I $ or to $ A \times \I $, with
   $ A $ an annulus;
   \item to each incidence of an edge $ e $ to a vertex $ v $ such that
   $ M_e $ is homeomorphic to $ T^2 \times \I $ we can associate an
   orientation reversing homeomorphism of a boundary component of $ M_e
   $ with a boundary component of $ M_v $ (note that in this case $ M_v
   $ must be Seifert fibered);
   \item to each incidence of an edge $ e $ to a vertex $ v $ such that
   $ M_e $ is homeomorphic to $ A \times \I $ we can associate an
   orientation reversing homeomorphism of a connected component of $
   (\partial A) \times \I $ with a subannulus of $ \partial M_v $ which
   is saturated when $ M_v $ is Seifert fibered and is a connected
   component of $ (\partial S_v) \times \I $ when $ M_v \cong S_v
   \times \I $;
   \item for all vertex manifolds $ M_v \cong S_v \times \I $ the
   valence of $ v $ is equal to the number of connected components of $
   (\partial S_v) \times \I $;
\end{itemize}
in order to have
$$ \frac{\bigcup_{v \in V} M_v \cup \bigcup_{e \in E} M_e}{\sim} \cong M. $$

The pair $ (M,\Gamma) $ is called a \emph{generalized graph manifold}. When a {\rm 3}-manifold $ M $ will be given a fixed generalized graph structure, we will abuse language and call $ M $ a generalized graph manifold.
A portion of the boundary of a vertex manifold $ M_v $ which is contained in $ \partial M $ will be called an \emph{exterior boundary}, and the union of all exterior boundaries of $ M_v $ will be denoted by $ \partial_M M_v $.

A generalized graph manifold will be called \emph{totally orientable} if all of its Seifert vertex manifolds are totally orientable as Seifert manifolds, and all $ S_v $'s are orientable surfaces.
\end{defi}

In \cite[Theorem 1.1]{schultens-graph} Schultens proves the following.

\begin{thm}[Schultens] \label{strirrsplitgengraph}
Let $ M $ be a totally orientable generalized graph manifold, and $ M = V \cup_\Sigma W $ a strongly irreducible Heegaard splitting. Then $ \Sigma $ can be isotoped so that:

\begin{enumerate}
   \item for each vertex manifold $ M_v $ the surface $ \Sigma \cap M_v
   $ is either horizontal, vertical, pseudohorizontal or pseudovertical;
   \item for each edge manifold $ M_e $ one of the following
   conditions holds:
   \begin{enumerate}
      \item[(i).] $ \Sigma \cap M_e $ is composed of incompressible annuli
      or is obtained from such a collection by {\rm 1}-surgery along an
      arc isotopic into $ \partial M_e $;
      \item[(ii).] $ M_e \cong T^2 \times \I $ and there are simple closed
      curves $ c, c' \subset T^2 \times \I $ such that $ c \cap c' = \{
      p \} $ and either $ V \cap (T^2 \times \I) $ or $ W \cap (T^2
      \times \I) $ is a collar of $ (c \times \{ 0 \}) \cup (\{ p \}
      \times \I) \cup (c' \times \{ 1 \}) $.   
   \end{enumerate}
\end{enumerate}
\end{thm}

\subsection{Manifolds with boundary}


Let $ M $ be an orientable Seifert fibered manifold with nonorientable base $ S $ and nonempty boundary, and fix a boundary component $ B $ of $ S $. Call $ g $ the genus of $ S $. Let $ \alpha_1, ..., \alpha_g \subset S $ be properly embedded arcs having boundary on $ B $ such that $ S \setminus (\alpha_1 \cup ... \cup \alpha_g) $ is a disk with holes, taking care that no arc passes through an exceptional point.

Let $ A_i $ be the annulus $ \pi^{-1}(\alpha_i) $ for $ i = 1, ..., g $. These have boundary on the torus $ T = \pi^{-1}(B) $; when we cut $ M $ along the $ A_i $'s we get a Seifert manifold $ M_0 $ with base a planar surface and the same exceptional fibers as $ M $: call $ A_i', A_i'' \subset \partial M_0 $ the annuli obtained when cutting along $ A_i $. For all $ i $ now connect $ A_i' $ to $ A_i'' $ with a copy of $ A_i \times \I $: as this retracts over $ A_i $, the result is the same manifold $ M $.

Now consider the graph $ \Gamma $ given by the wedge of $ g $ copies of $ \esseuno $, with vertex set $ V = \{ v_0 \} $ and edge set $ E = \{ e_1, ..., e_g \} $. We associate $ M_0 $ to $ v_0 $ and $ A_i \times \I $ to $ e_i $ for all $ i $; then we realize the incidence of $ e_i $ on $ v_0 $ by gluing $ A_i' $ to $ A_i \times \{ 0 \} $ via the identity map and $ A_i \times \{ 1 \} $ to $ A_i'' $ via the orientation reversing map that identifies $ A_i' $ with $ A_i'' $.

We have thus given $ (M, \Gamma) $ a totally orientable generalized graph structure. Note that the very same construction works when $ M $ is totally orientable, and the graph $ \Gamma $ is nontrivial in this case.
\begin{thm} \label{strirrboundary}
Let $ M $ be a Seifert fibered manifold with nonorientable base $ S $ and nonempty boundary. Then every strongly irreducible Heegaard splitting of $ M $ is vertical.
\end{thm}
\begin{proof}
Let $ \Sigma $ be the Heegaard surface, and equip $ M $ with the totally orientable generalized graph structure described above, denoting again with $ M_0 $ the vertex manifold. We apply Theorem \ref{strirrsplitgengraph}, and find an isotopy of $ \Sigma $ so that the following conditions are satisfied:
\begin{enumerate}
   \item $ \Sigma \cap M_0 $ is either horizontal, vertical,
   pseudohorizontal or pseudovertical. Note that any horizontal or
   pseudohorizontal surface must be horizontal near $ \partial_M M_0
   $, in particular $ \Sigma $ would meet $ \partial M $, which is
   clearly absurd. Hence $ \Sigma \cap M_0 $ is either vertical or
   pseudovertical.
   \item For all $ i $ the intersection $ \Sigma_i = \Sigma \cap (A_i
   \times \I) $ either consists of incompressible annuli or is obtained
   from such a collection by {\rm 1}-surgery along one single arc
   isotopic into $ \partial (A_i \times \I) $. For the same reason we
   excluded (pseudo)horizontal components before, we see that each
   annulus in the collection has boundary in $ A_i \times \{ 0, 1 \} $.
\end{enumerate}
Note that possibility $ (ii) $ of Theorem \ref{strirrsplitgengraph} cannot occur because all edge manifolds are homeomorphic to $ (annulus) \times \I $.

Now fix an edge manifold $ A_i \times \I $. Note that a closed curve contained in $ A_i \times \{ 0, 1 \} $ can be extended locally into $ M_0 $ to a vertical surface. Hence, if we have no added arcs, then every incompressible and boundary parallel annulus in $ A_i \times \I $ can be pushed through $ A_i \times \{ 0, 1 \} $ using an innermost argument, and the (pseudo)verticality of $ \Sigma \cap M_0 $ is preserved. Hence, if we have no added arcs, we can assume that every connected component of $ \Sigma_i $ is an annulus which is vertical in the fibration $ A_i \times \I \rightarrow A_i $.

If we have an added arc $ \alpha $ in $ A_i \times \I $, note first that it must have endpoints on different components of $ \Sigma_i $, otherwise $ \Sigma $ would be stabilized. If $ \alpha $ connects two annuli which are vertical in the fibration $ A_i \times \I \rightarrow A_i $ we may push it through $ A_i \times \{ 0 \} $ to a horizontal arc in $ M_0 $. If $ \Sigma \cap M_0 $ was vertical, it is now pseudovertical; if it was pseudovertical this condition is preserved, as we may keep $ \alpha $ in a neighbourhood of $ A_i' $. So we are reduced to the case when we have no arcs in $ A_i \times \I $. The same holds if $ \alpha $ connects a vertical annulus to a boundary parallel one.

\begin{figure}[!htb]
\begin{center}
\input{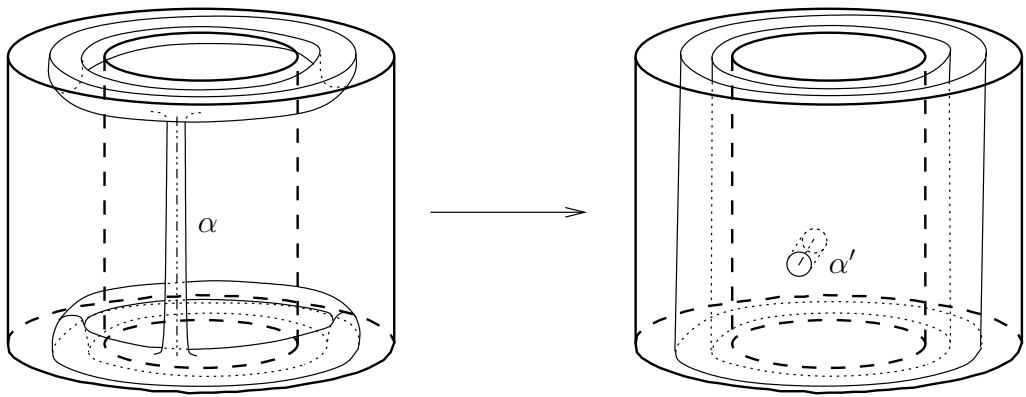}
\caption{A move analogous to the swap described in Lemma \ref{lemswap}.} \label{arcsinedges}
\end{center}
\end{figure}

On the other hand, if $ \alpha $ connects two boundary parallel annuli, we have two cases:
\begin{enumerate}
   \item the two annuli have boundary on the same connected component
   of $ A_i \times \{ 0, 1 \} $. We can then push $ \alpha $ through $
   A_i \times \{ 0, 1 \} $ as before.
   \item the two annuli have boundary on different components of $ A_i
   \times \{ 0, 1 \} $. Then we can isotope the corresponding component
   of $ \Sigma_i $ until it is obtained by adding a ``dual'' horizontal
   arc $ \alpha' $ to two vertical annuli (see Figure \ref
   {arcsinedges}). Again, we reduce ourselves to the case when we have
   no arcs in $ A_i \times \I $.
\end{enumerate}

To sum up, we have found an isotopy that puts $ \Sigma \cap M_0 $ in (pseudo)vertical position, and the $ \Sigma_i $'s all consist of vertical annuli. When we retract the edge manifolds $ A_i \times \I $ onto $ A_i $ to recover $ M $, the effect is to retract every connected component of $ \Sigma_i $ to a fiber of $ M $. That is, $ \Sigma $ is a vertical or pseudovertical surface. As it is connected, in the former case it is a torus, i.e. $ M $ admits a genus $ 1 $ Heegaard splitting and hence $ M $ is a lens space, whose splittings are all vertical by \cite{bonahonotal-heegaardlens}. In the latter case, Proposition \ref{equivvert} implies that the splitting is vertical.
\end{proof}

\subsection{Weakly reducible splittings}

Now $ M $ will denote an orientable Seifert manifold with nonorientable base, perhaps with boundary.

\begin{lem}[Jaco] \label{incompressible}
Let $ F $ be a two-sided incompressible surface in $ M $. Then either $ F $ is a vertical torus or annulus, or it splits $ M $ into two parts both homeomorphic to some twisted product $ B \tilde{\times} \I $, or is boundary parallel.
\end{lem}
\begin{proof}
According to \cite[Theorem VI.34]{jaco-3mfld}, we only need to show that there is no fibration $ M \rightarrow \esseuno $ which admits $ F $ as a fiber. But such a fibration would induce an orientation on the Seifert fibration we have on $ M $, and this is a contradiction because a ``positive'' vector along the fiber would induce an orientation on the base $ S $.
\end{proof}

We now have enough tools to characterize irreducible but weakly reducible Heegaard splittings of $ M $. The proof of the following Theorem mimics that of \cite[Theorem 2.6]{moriahschultens-seifert}.

\begin{thm} \label{irrweakredsplit}
Let $ M $ be a Seifert manifold with nonorientable base and $ \Sigma $ an irreducible Heegaard surface for $ M $. If $ \Sigma $ is weakly reducible, then it is vertical.
\end{thm}
\begin{proof}
Thanks to Remark \ref{onlylarge}, we can always find a curve $ \gamma \subset S $ that avoids all exceptional points in $ S $ and its preimage under the projection map is an incompressible torus $ T $.

As $ \Sigma $ is weakly reducible, there are disjoint families of essential compressing disks $ \Delta_i \subset W_i $ for $ i = 1, 2 $ such that the surface $ \Sigma^* $ obtained compressing $ \Sigma $ along $ \Delta_1 \cup \Delta_2 $ is incompressible, and each connected component of $ M \setminus \Sigma^* $ inherits a strongly irreducible Heegaard splitting from $ M $ (recall Theorem \ref{weakreduction}). We choose the $ \Delta_i $'s with the additional constraint that $ (genus(\Sigma^*), | \Sigma^* \cap T |) $ is minimal according to lexicographical order.

Let $ \Sigma^*_1 $ be a connected component of $ \Sigma^* $. Thanks to Lemma \ref{incompressible}, we have two possibilities for $ \Sigma^*_1 $:
\begin{enumerate}
   \item $ \Sigma^*_1 $ is a horizontal surface, and separates $ M $
   into two twisted $ \I $-bundles over a surface $ F $. In this case
   all components of $ \Sigma^* $ are parallel.
   
   Let $ c $ be a connected component of $ \Sigma^*_1 \cap T $; as both
   surfaces are incompressible and the number of intersections of $ T $
   and $ \Sigma^* $ was chosen to be minimal, $ c $ is essential in
   both. Let $ N_i $ ($ i = 1, 2 $) be the connected components of $ M
   \setminus \Sigma^*_1 $ such that $ c \subset \partial N_1 \cap \partial
   N_2 $. Then either $ N_i \cong \Sigma^*_1 \times \I $ or $ N_i \cong
   \fti $, with $ F $ a nonorientable surface.
   
   If at least one of the $ N_i $'s is a product, then $ \Sigma^* $ has
   at least two connected components and it is easy to see that $ \Sigma $ is
   reducible.
   
   If both $ N_i $'s are twisted products, then $ \Sigma^*_1 $ is the
   only connected component of $ \Sigma^* $, and $ M $ is built by
   gluing together two copies of $ \fti $ along some finite-order
   homeomorphism between their boundaries. We have exactly one vertical arc $
   \alpha_i $ embedded in each of the $ N_i $'s which induces a Heegaard
   splitting on it: we choose $ \alpha_i $ to lie in the vertical annulus
   $ A_i $ which has $ c $ as boundary and which is contained in $
   T \cap N_i $, and we make sure that $ (\alpha_1 \cap c) \cap
   (\alpha_2 \cap c) = \emptyset $. When we reconstruct $ \Sigma $ from
   $ \Sigma^* $, the $ A_i $'s become disjoint disks $ D_i $'s lying in
   different sides of $ \Sigma $.

   \begin{figure}[!htb]
   \begin{center}
   \input{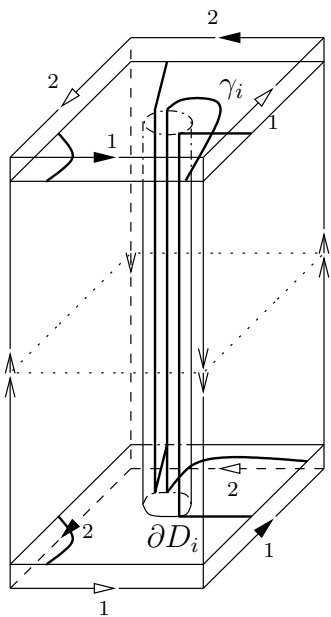}
   \caption{The curves $ \gamma_i $ and $ \partial D_i $ meet in a single
   point.} \label{weakridrid2}
   \end{center}
   \end{figure}
   
   Note that as $ T $ is two-sided in $ M $, the $ D_i $'s can be pushed
   away from $ T $ so that $ D_1 \cap T = D_2 \cap T = \emptyset $.
   In addition, Figure \ref{weakridrid2} shows that $ \partial D_i $ are
   nonseparating in $ \Sigma^* $, hence they are essential (when $ F $ is
	not the Klein bottle one sees that the same choice of curves works by
	constructing $ F $ as $ K^2 \# F' $, thus removing from both
	$ K^2 \tilde{\times} \I $ and $ F' \tilde{\times} \I $ a neighbourhood
	of a fiber, and then pasting the two manifolds via a fiber preserving
	homeomorphism; $ F $ cannot be a projective plane, or the ramified
	covering map $ F \rightarrow S $ would be a genuine homeomorphism, and
	its inverse would be a section of the fibration, hence we would have
	$ M \cong \pidue \tilde{\times} \esseuno $ against our assumption).
   
   Now compress $ \Sigma $ along $ D_1 \cup D_2 $. This gives a new
   surface which might or not be compressible. If it is, we may
   compress it further to get an incompressible surface whose genus is
   smaller than $ \Sigma^* $'s. If it is incompressible, we have an
   incompressible surface with genus equal to $ \Sigma^* $'s that
   intersects $ T $ fewer times that $ \Sigma^* $ does. In any case, we
   have a contradiction thanks to the choice of $ \Delta $. Hence no
   component of $ \Sigma^* $ is a horizontal surface.
   
   \item All components of $ \Sigma^* $ are then vertical tori. Thanks
   to Theorem \ref{weakreduction}, $ \Sigma $ is obtained by
   amalgamation of Heegaard splittings of the connected components of
   $ M \setminus \Sigma^* $. As these manifolds have boundary, all
   Heegaard splittings induced there are vertical thanks to Theorem
   \ref{strirrboundary}. Hence $ \Sigma $ gives a vertical Heegaard
   splitting thanks to Proposition \ref{amalgvert}.
\end{enumerate}
The proof is complete.
\end{proof}

\subsection{Horizontal splittings}
It is straightforward from the definition that every Seifert manifold admits vertical Heegaard splittings. We shall now describe, following \cite{moriahschultens-seifert}, another construction that may give rise to a Heegaard splitting.

Let $ M^* $ be a Seifert manifold with nonorientable base and one boundary component, and let $ S $ be a horizontal two-sided surface (which always exists as $ \partial M^* \neq \emptyset $). Thanks to Lemma \ref{incompressible}, $ S $ splits $ M^* $ into two twisted $ \I $-bundles $ W_1^* $, $ W_2^* $ over surfaces $ F_1 \cong F_2 $: each of the $ W_i^* $'s is a handlebody. The annuli $ A_i = \partial F_i \tilde{\times} \I $, for $ i = 1, 2 $, are glued together to build up $ \partial M^* $.

Fix a Dehn filling of $ \partial M^* $, and suppose we obtain a Seifert manifold $ M $. Set $ W_1 = W_1^* $ and $ W_2 = W_2^* \cup V $, where $ V $ is the solid torus in the Dehn filling, and let $ \Sigma = \partial W_1 = \partial W_2 $. Note that $ \Sigma $ embeds in $ M $. In order for $ \Sigma $ to be a splitting surface for $ M $, $ W_2 $ must be a handlebody. As in the totally orientable case, the right condition for this to happen is to choose the Dehn filling to glue $ \partial F_2 $ to a generator of $ \pi_1(V) $. Note that this is equivalent to say that the gluing must have a matrix
$ \left(
\begin{array}{cc}
a & b \\
c & d
\end{array}
\right) \in GL(2,\mathbb{Z}) $
with $ a = 1 $ in the framing determined by $ (\partial F_1, $ \emph{regular fiber}$ ) $.

\begin{defi} \label{horizontal}
Let $ M $ be a Seifert fibered manifold with nonorientable base, and $ f $ a fiber. Let $ S $ be a horizontal surface in $ M^* = M \setminus \eta(f) $. If $ M $ is constructed from $ M^* $ via a $ (1,n) $ Dehn filling, then the Heegaard splitting constructed above is called a \emph{horizontal} splitting corresponding to the fiber $ f $.
\end{defi}

\subsection{Closed Seifert manifolds}

For the rest of the paper $ M $ will be a Seifert manifold with nonorientable base and no boundary. This section is devoted to the study of Heegaard splittings of such manifolds, and it follows very closely \cite[Theorem 0.1]{moriahschultens-seifert}.

Let $ x_1, ..., x_m \in S $ be the exceptional points on $ S $; if there are no exceptional fibers fix a regular point and call it $ x_1 $. Let $ D_1, ..., D_m $ be pairwise disjoint disk neighbourhoods of $ x_1, ..., x_m $, and set $ S^* = S \setminus (D_1 \cup ... \cup D_m) $. Fix $ p \in S^* $ and choose simple closed curves $ \alpha_1, ..., \alpha_g $ based in $ p $ such that the surface obtained by cutting $ S^* $ along $ \bigcup_{j = 1}^g \alpha_j $ is a disk with holes. In addition choose simple closed curves $ c_1, ..., c_m $ also based in $ p $ such that $ c_i $ cuts off a regular neighbourhood of $ \partial D_i $ in $ S^* $, for $ i = 1, ..., m $.

As we removed a neighbourhood of the exceptional points (or of a regular point in which we concentrate the Euler number if there are no exceptional fibers), the surface $ S^* $ embeds in $ M $, and the manifold $ M^* = \pi^{-1}(S^*) $ is Seifert with base $ S^* $. Consider the preimages $ A_j = \pi^{-1}(\alpha_j) $ for $ j = 1, ..., g $ and $ C_i = \pi^{-1}(c_i) $ for $ i = 1, ..., m $. The former is a collection of vertical Klein bottles and the latter is a collection of vertical tori, and each element embeds in $ M $. Now cut $ M $ along $ \bigcup_{j = 1}^g A_j \cup \bigcup_{i = 1}^m C_i $. We obtain a collection $ V_0, ..., V_m $ of solid tori, where the exceptional fiber $ f_i $ is the core of the solid torus $ V_i $ for $ i = 1, ..., m $. There is one copy of the regular fiber $ \pi^{-1}(p) $ on each boundary $ \partial V_i $ for $ i = 1, ..., m $, and $ 2g $ copies on $ \partial V_0 $, and all of them are longitudes of the corresponding solid torus.

Call $ \tilde{A_j} $ and $ \tilde{C_i} $ the annuli obtained from $ A_j $ and $ C_i $ respectively when cutting $ \partial V_i $, for $ i = 0, ..., m $ along the copies of $ \pi^{-1}(p) $.

The following two Lemmas are proved verbatim by the proofs of \cite[Lemma 1.5]{moriahschultens-seifert} and \cite[Lemma 1.6]{moriahschultens-seifert}.

\begin{lem} \label{noiso1}
If there is no isotopy pushing a fiber into the Heegaard surface $ \Sigma $, then there is an isotopy such that all transversal intersections $ \Sigma \cap \tilde{A_j} $ and $ \Sigma \cap \tilde{C_i} $ contain at least one essential arc and no inessential arcs.
\end{lem}

\begin{lem} \label{noiso2}
In the same hypotesis of Lemma \ref{noiso1}, the simple closed curves composing $ \Sigma \cap \partial V_i $ are meridians of $ \partial V_i $, for $ i = 0, ..., m $.
\end{lem}

The following Proposition characterizes strongly irreducible Heegaard splittings. The same proof of \cite[Proposition 1.3]{moriahschultens-seifert} works here.

\begin{prop} \label{weakredhasiso}
Let $ \Sigma $ be a Heegaard splitting of $ M $. Then either there is an isotopy pushing a fiber into $ \Sigma $, or $ \Sigma $ is weakly reducible.
\end{prop}

We can finally prove the main Theorem of this paper. The proof mimics that of \cite[Theorem 0.1]{moriahschultens-seifert}.

\begin{thm} \label{main}
Let $ \Sigma $ be an irreducible Heegaard splitting of a closed Seifert manifold with nonorientable base space. Then it is either vertical or horizontal.
\end{thm}
\begin{proof}
Thanks to Theorem \ref{irrweakredsplit}, we can assume that $ \Sigma $ is a strongly irreducible Heegaard surface.

Thanks to Proposition \ref{weakredhasiso}, we may perform an isotopy that brings a fiber $ f $ into $ \Sigma $. Remove an open neighbourhood $ \eta(f) $ of $ f $ from $ M $ to obtain a manifold $ M^* $, and set $ \Sigma^* = \Sigma \cap M^* $. We now distinguish two cases.

\underline{Case 1}: $ \Sigma^* $ is incompressible in $ M^* $. In this case, Lemma \ref{incompressible} implies that it is either a vertical annulus, or it splits $ M^* $ into two twisted $ \I $-bundles over a nonorientable surface $ B $; $ \Sigma $ is constructed by gluing the annulus $ \Sigma \cap \overline{\eta(f)} $ to $ \Sigma^* $.

Hence, in the former case $ \Sigma $ is a vertical torus and so $ M $ is a lens space and admits a totally orientable Seifert structure with base $ \essedue $. As the result is already known in this case (\cite{bonahonotal-heegaardlens}, \cite{moriahschultens-seifert}), we may assume that this does not happen.
   
In the latter case, the handlebody $ W_2 $ may be constructed by gluing a solid torus to one component of $ M^* \setminus \Sigma^* $. As the result yields a handlebody, this gluing must be a $ (1,n) $-Dehn filling. Hence, by Definition \ref{horizontal}, $ \Sigma $ gives a horizontal Heegaard splitting.
   
\underline{Case 2}: $ \Sigma^* $ is compressible in $ M^* $. In this case, choose a collection $ \Delta \subset M^* $ of disjoint compressing disks for $ \Sigma^* $ with the constraint that it minimizes the number of intersections with $ \Sigma^* $: it must be contained entirely in one side of $ \Sigma^* $, say $ \Delta \subset W_2 $, otherwise $ \Sigma $ would be weakly reducible (recall Lemma \ref{weaklyred}).

Let $ \Sigma^{**} $ be the surface obtained compressing $ \Sigma^* $ along $ \Delta $, and call $ \Sigma^{**}_0 $ one of its connected components with nonempty boundary. Then Lemma \ref{incompressible} implies that $ \Sigma^{**}_0 $ can either be an annulus, or it can split $ M^* $ into two twisted $ \I $-bundles (and in either case it is the only component of $ \Sigma^{**} $ with nonempty boundary).

Suppose $ \Sigma^{**}_0 $ is an annulus. Then there are two possibilities:
\begin{itemize}
   \item either $ \Sigma^{**}_0 $ is boundary parallel, and it cuts off a trivially fibered solid torus, or
   \item $ \Sigma^{**}_0 \cup (\Sigma \cap \overline{\eta(f)}) \subset W_2 $ is a solid torus which is nontrivially fibered.
\end{itemize}
Note that $ \Sigma^{**}_0 \cup (\Sigma \cap \overline{\eta(f)}) $ cannot be an incompressible torus, because it can be isotoped entirely into $ W_2 $.

In both cases we may perform a small isotopy so to push $ f $ into $ W_2 $, so that it becomes a core. We may now remove a neighbourhood of $ f $ from $ M $ to obtain a manifold homeomorphic to $ M^* $ that admits $ \Sigma $ as a splitting surface. This splitting is vertical by Theorems \ref{strirrboundary} and \ref{irrweakredsplit}, and it is straightforward from Definition \ref{defvertical} that the given splitting is also vertical.

\linea

We claim that $ \Sigma^{**}_0 $ cannot split $ M^* $ into two twisted $ \I $-bundles $ W_1^{**} $, $ W_2^{**} $ (which are handlebodies, as the surface they fiber on has boundary), or $ \Sigma $ would be reducible. In fact, $ \Sigma $ induces a Heegaard splitting on $ W_1^{**} $ as follows (recall that $ \Delta \subset W_2 $, and that we get $ W_1 $ from $ W_1^{**} $ by drilling out tunnels along the cocores of $ \Delta $): first push $ \Sigma \cap \overline{\eta(f)} $ through $ \partial M^* $, then glue the {\rm 1}-handles that compose $ W_2 \cap W_1^{**} $ to a collar of $ \partial W_1^{**} $ in $ W_1^{**} $ to obtain a compression body $ W_2' $, and set $ W_1' = \overline{W_1^{**} \setminus W_2'} $. As $ \Delta $ was entirely contained in $ W_2 $, $ (W_1',W_2') $ is a Heegaard splitting of $ W_1^{**} $; this splitting cannot be trivial, because $ \Delta \neq \emptyset $. Hence it is stabilized by Theorem \ref{handlebody}; as $ \Sigma $ is constructed by gluing two-dimensional {\rm 1}-handles to $ \Sigma^{**}_0 $, the given splitting is also stabilized, contradicting our hypotesis.
\end{proof}

\nocite{cassongordon-reducingheegaard, schultens-weakly, boileauotal-heegaardt3}
\bibliographystyle{siam}
\bibliography{biblio}

\end{document}